\documentclass[12pt,a4paper]{article}
\usepackage{enumerate}

\usepackage{amsmath,amssymb,xspace,amsthm}
\newtheorem{theorem}{Theorem}%[section]
\newtheorem{lemma}[theorem]{Lemma}%[section]
\newtheorem{corollary}[theorem]{Corollary}%[section]
\numberwithin{equation}{section}

\newcommand{\C}{\ensuremath{\Bbbk}\xspace}

\newcommand{\h}{\ensuremath{\mathfrak{H}}}

\newcommand{\Z}{\ensuremath{\mathbb{Z}}\xspace}

\newcommand{\V}{\ensuremath{\mathcal{V}}\xspace}

\renewcommand{\phi}{\varphi}

\renewcommand{\geq}{\geqslant}

\begin{document}
\title{Classification of simple weight Virasoro modules with a finite-dimensional weight space}
\author{Volodymyr Mazorchuk and Kaiming Zhao}
\date{}
\maketitle



\begin{abstract}
We show that the support of a simple weight module over the Virasoro
 algebra, which has an infinite-dimensional weight space, coincides
with the weight lattice and that all non-trivial weight spaces of
such module are infinite dimensional. As a corollary we obtain
that every simple weight module over the Virasoro algebra, having
a non-trivial finite-dimensional weight space, is a Harish-Chandra
module (and hence is either a simple highest or lowest weight
module or a simple module from the intermediate series). This
implies positive answers to two conjectures about simple pointed
and simple mixed modules over the Virasoro algebra.
\end{abstract}

\section{Description of the results}\label{s1}

The {\em Virasoro algebra} $\V$ over an algebraically closed
field, $\C$, of characteristic zero has a basis, consisting of a
central element, $c$, and elements $e_i$, $i\in\Z$, with the Lie
bracket defined for the basis elements as follows:
\begin{displaymath}
[e_i,e_j]=(j-i)e_{i+j}+\delta_{i,-j}\frac{i^3-i}{12}c.
\end{displaymath}
The linear span of $c$ and $e_0$ is called the {\em Cartan subalgebra}
$\h$ of $\V$ and an $\h$-diagonalizable $\V$-module is usually called a
{\em weight module}. If, additionally, all weight spaces of a weight $\V$-module are
finite-dimensional, the module is called a {\em Harish-Chandra module}, see for
example \cite{M}. All simple Harish-Chandra modules were classified in
\cite{MP1,MP2,M} and are exhausted by simple highest weight modules, simple
lowest weight modules and simple modules from the so-called {\em intermediate
series} (see e.g. \cite{M} for definitions).

If $M$ is a simple weight $\V$-module, then $c$ acts on $M$ by a scalar, called
the {\em central charge} of $M$. Furthermore, $M$ can be written as a direct
sum of its {\em weight spaces}, $M=\oplus_{\lambda\in\C}M_{\lambda}$,
where $M_{\lambda}$ is the set of all elements of $M$ on which $e_0$ acts
as the multiplication with $\lambda$. The set of all $\lambda$ for which
$M_{\lambda}\neq 0$ is called the {\em support} of $M$ and is denoted by
$\mathrm{supp}(M)$. Obviously, if $M$ is a simple weight $\V$-module, then there
exists $\lambda\in\C$ such that $\mathrm{supp}(M)\subset\lambda+\Z$. A simple
weight $\V$-module, $M$, is called {\em pointed} provided that there exists
$\lambda\in\C$ such that $\dim M_{\lambda}=1$ (for example from the classification
of simple Harish-Chandra modules it follows that they all are pointed). The following
question was formulated in \cite[Problem~3.3]{X}:
\vspace{1mm}

{\bf Question:} {\em Is any simple pointed $\V$-module a Harish-Chandra module?}
\vspace{1mm}

A simple weight $\V$-module, $M$, is called {\em mixed} provided that there
exist $\lambda\in\C$ and $k\in \Z$ such that $\dim M_{\lambda}=\infty$ and
$\dim M_{\lambda+k}<\infty$. The following
conjecture, a positive answer to which implies a positive answer to the
question above, was formulated in \cite[Conjecture~1]{Ma1}:
\vspace{1mm}

{\bf Conjecture:} {\em There are no simple mixed $\V$-modules.}
\vspace{1mm}

In the present paper we prove the following result, which implies positive
answers to both the Question and the Conjecture above:

\begin{theorem}\label{tmain}
Let $M$ be a simple weight $\V$-module. Assume that there exists
$\lambda\in \C$ such that $\dim M_{\lambda}=\infty$. Then
$\mathrm{supp}(M)=\lambda+\Z$ and for every $k\in\Z$ we have
$\dim M_{\lambda+k}=\infty$.
\end{theorem}

Apart from the positive answers to the Question and the Conjecture
above, Theorem~\ref{tmain} also implies the following
classification of all simple weight $\V$-modules which admit a
non-trivial finite-dimensional weight space:

\begin{corollary}\label{cmain}
Let $M$ be a simple weight $\V$-module. Assume that there exists
$\lambda\in \C$ such that $0<\dim M_{\lambda}<\infty$. Then $M$ is
a Harish-Chandra module. Consequently, $M$ is either a simple
highest or lowest weight module or a simple module from the
intermediate series.
\end{corollary}

The paper is organized as follow: Theorem~\ref{tmain} is proved in Section~\ref{s2}
and in Section~\ref{s3} we discuss the corollaries from this theorem.
\vspace{0.2cm}

\noindent
{\bf Acknowledgments.} The research was done during the visit of the first author to
Beijing Normal University in May 2005. The hospitality and
financial support of Beijing Normal University are gratefully
acknowledged. For the first author the research was partially
supported by the Royal Swedish Academy of Sciences, the Swedish
Research Council, and STINT (the Swedish Foundation for International
Cooperation in research and Higher Education). The second
author is supported by  NSERC, and NSF of China (Grant 10371120
and 10431040).

\section{Proof of Theorem~\ref{tmain}}\label{s2}

Our strategy to prove Theorem~\ref{tmain} is the following: first we show in
Lemma~\ref{l1} that a simple weight $\V$-module with an infinite-dimensional
weight space can have at most one weight space of finite dimension in the weight 
lattice. Then in Lemma~\ref{l3} we show that this finite-dimensional weight must
belong to $\{-1,0,1\}$. These three cases are excluded in Lemma~\ref{l4} and
Lemma~\ref{l5} by a case-by-case analysis. The key point of our proof is
Lemma~\ref{l2}, which is an easy technical statement claiming that some
special element from the universal enveloping algebra $U(\V)$ annihilates
certain elements of the module. The statement is very easy to prove by a
direct computation, however, the main and perhaps most non-trivial idea of
the paper is that there should exist an element in the universal enveloping
algebra $U(\V)$, which satisfies the conclusion of Lemma~\ref{l2}.

\vspace{2mm}
Let $M$ be a simple weight $\V$-module. We start with
the following obvious observation: \vspace{2mm}

\noindent
{\bf Principal Observation:} {\em Assume that there exists $\mu\in\C$ and a
non-zero element, $v\in M_{\mu}$, such that $e_1v=e_2v=0$ or $e_{-1}v=e_{-2}v=0$. Then $M$
is a Harish-Chandra module.}
\vspace{2mm}

\begin{proof}
Indeed, under these conditions $v$ is either a highest or a lowest weight vector
and hence $M$ is either a highest or a lowest weight module. Hence $M$ is a
Harish-Chandra module (see e.g. \cite{M}).
\end{proof}

Assume now that $M$ is a simple weight $\V$-module and that there exists
$\lambda\in \C$ such that $\dim M_{\lambda}=\infty$.

\begin{lemma}\label{l1}
There exists at most one $i\in\Z$ such that $\dim M_{\lambda+i}<\infty$.
\end{lemma}

\begin{proof}
Assume that $\dim M_{\lambda+i}<\infty$ and $\dim M_{\lambda+j}<\infty$ for
different $i,j\in\Z$. Without loss of generality we may assume
$i=1$ and $j>1$. Let $V$ denote the intersection of the kernels of the linear
maps $e_1:M_{\lambda}\to M_{\lambda+1}$ and $e_j:M_{\lambda}\to M_{\lambda+j}$.
Since $\dim M_{\lambda}=\infty$, $\dim M_{\lambda+1}<\infty$ and
$\dim M_{\lambda+j}<\infty$, we have $\dim V=\infty$. Since
$[e_1,e_k]=(k-1)e_{k+1}\neq 0$ for all $k>1$, and $j>1$, inductively we get
\begin{equation}\label{eq1}
e_k V=0\quad\text{ for all  }\quad k=1,j,j+1,j+2,\dots.
\end{equation}

(Here we cannot directly use the well known \cite[Lemma~1.5]{M} to
deduce that $M$ is a highest weight module).  If there would exist
$0\neq v\in V$ such that $e_2v=0$, then $e_1v=e_2v=0$ and $M$
would be a Harish-Chandra module by the Principal Observation. A
contradiction. Hence $e_2v\neq 0$ for all $v\in V$. In particular,
$\dim e_2 V=\infty$. Since $\dim M_{\lambda+1}<\infty$, there
exists $0\neq w\in e_2 V$ such that $e_{-1}w=0$. Let $w=e_2u$ for
some $u\in V$. For all $k\geq j$, using \eqref{eq1} we have
\begin{displaymath}
e_kw=e_ke_2u=e_2e_ku+(2-k)e_{k+2}u=0+0=0.
\end{displaymath}
Hence $e_kw=0$ for all $k=-1,j,j+1,j+2,\dots$. Since
$[e_{-1},e_l]=(l+1)e_{l-1}\neq 0$ for all $l>1$, inductively we get
$e_kw=0$ for all $k=1,2,\dots$. Hence $M$ is a Harish-Chandra module by
the Principal Observation. A contradiction. The statement follows.
\end{proof}

Because of Lemma~\ref{l1} we can now fix the following notation until
the end of this section: $M$ is a simple weight $\V$-module, $\mu\in\C$ is such that
$\dim M_{\mu}<\infty$ and $\dim M_{\mu+i}=\infty$ for every $i\in\Z\setminus\{0\}$.

\begin{lemma}\label{l2}
Let $0\neq v\in M_{\mu-1}$ be such that $e_1v=0$. Then
\begin{displaymath}
(e_1^3-6e_2e_1+6e_3)e_2v=0.
\end{displaymath}
\end{lemma}

\begin{proof}
\begin{multline*}
(e_1^3-6e_2e_1+6e_3)e_2v=(e_1^3e_2-6e_2e_1e_2+6e_3e_2)v=\\
(e_2e_1^3+3e_3e_1^2+6e_4e_1+6e_5-6e_2^2e_1-6e_3e_2-6e_5+6e_3e_2)v=\\
(e_2e_1^3+3e_3e_1^2+6e_4e_1-6e_2^2e_1)v=[\text{using }e_1v=0]=0.
\end{multline*}
\end{proof}

\begin{lemma}\label{l3}
$\mu\in\{-1,0,1\}$.
\end{lemma}

\begin{proof}
Let $V$ denote the kernel of $e_1:M_{\mu-1}\to M_{\mu}$. Since
$\dim M_{\mu-1}=\infty$ and $\dim M_{\mu}<\infty$ we have $\dim
V=\infty$. For any $v\in V$ consider the element $e_2v$. By the
Principal Observation, $e_2v=0$ would imply that $M$ is a
Harish-Chandra module, a contradiction. Hence $e_2v\neq 0$, in
particular, $\dim e_2V=\infty$. This implies that there exists
$w\in e_2V$ such that $w\neq 0$ and $e_{-1}w=0$. From
Lemma~\ref{l2} we have $(e_1^3-6e_2e_1+6e_3)w=0$, in particular,
$e_{-1}^3(e_1^3-6e_2e_1+6e_3)w=0$. However, by a direct
calculation we obtain
\begin{displaymath}
e_{-1}^3(e_1^3-6e_2e_1+6e_3)=48e_0^3-144e_0^2+96e_0\mod U(\V)e_{-1}.
\end{displaymath}
This implies $(48e_0^3-144e_0^2+96e_0)w=0$. But $w\in M_{\mu+1}$, which
implies $e_0w=(\mu+1)w$, and hence
$(\mu+1)^3-3(\mu+1)^2+2(\mu+1)=0$, that is $\mu\in\{-1,0,1\}$.
\end{proof}

\begin{lemma}\label{l4}
$\mu\in\{-1,1\}$ is not possible.
\end{lemma}

\begin{proof}
We show that $\mu=1$ is not possible and for $\mu=-1$ the
statement will follow by applying the canonical involution on
$\V$. Assume $\mu=1$ and denote by $V$ the infinite-dimensional
kernel of the linear map $e_1:M_0\to M_1$. For  $v\in V$, using
$e_1v=e_0v=0$ we have
\begin{equation}\label{eq2}
e_1e_{-1}v=e_{-1}e_{1}v-2e_0v=0+0=0.
\end{equation}
Hence if $e_{-1}V$ would be infinite-dimensional, there would
exist $0\neq w\in e_{-1}V$ such that $e_1w=0$ (by \eqref{eq2}) and
$e_2w=0$ (since $\dim V_1<\infty$). The Principal Observation then
would imply that $M$ is a Harish-Chandra module, a contradiction.
Hence $\dim e_{-1}V<\infty$. This means that the kernel $W$ of the
linear map $e_{-1}:V\to M_{-1}$ is infinite-dimensional. For every
$x\in W$ we have
\begin{equation}\label{eq3}
e_1e_{-2}x=e_{-2}e_1x-3e_{-1}x=[\text{using }e_{-1}x=e_1x=0]=0+0=0.
\end{equation}
If there would exist $0\neq x\in W$ such that $e_{-2}x=0$, then we
would have $e_{-2}x=e_{-1}x=0$ and  the Principal Observation
would imply that $M$ is a Harish-Chandra module, a contradiction.
Thus $\dim e_{-2}W=\infty$. Let $H$ denote the kernel of the
linear map $e_3:e_{-2}W\to M_1$. Since $\dim e_{-2}W=\infty$ and
$\dim M_1<\infty$, we have $\dim H=\infty$. For every $y\in H$ we
also have $e_1y=0$ by \eqref{eq3}, implying by induction that $e_k
H=0$ for all $k=1,3,4,\dots$.

If $e_2 h=0$ for some $0\neq h\in H$ then the Principal
Observation implies that $M$ is a Harish-Chandra module, a
contradiction. Hence $\dim e_2H=\infty$. For every $h\in H$ and
$k\geq 3$ we have
\begin{displaymath}
e_ke_2h=e_2e_kh+(2-k)e_{k+2}h=[\text{using } e_ih=0\text{ for }i\geq 3]=0+0=0.
\end{displaymath}
Hence $e_ke_2H=0$ for all $k\geq 3$. Let, finally, $K$ denote the
infinite-dimensional kernel of the linear map $e_1:e_2H\to M_1$.
If $e_2z=0$ for some $0\neq z\in K$ then the Principal Observation
implies that $M$ is a Harish-Chandra module, a contradiction.
Hence $\dim e_2z\neq 0$ for all $z\in K$. For every $z\in K$ and
$k\geq 3$ we have
\begin{displaymath}
e_ke_2z=e_2e_kz+(2-k)e_{k+2}z=[\text{using } e_iz=0\text{ for }i\geq 3]=0+0=0.
\end{displaymath}
Hence $e_ke_2K=0$ for all $k\geq 3$. At the same time, since
$\dim e_2K=\infty$ and $\dim M_1<\infty$, we can find some $0\neq t\in e_2K$
such that $e_{-1}t=0$. By induction we get $e_it=0$ for all $i>0$ and thus
$M$ is a Harish-Chandra module by the Principal Observation. This last
contradiction completes the proof.
\end{proof}

Now the proof of Theorem~\ref{tmain} follows from the following lemma:

\begin{lemma}\label{l5}
$\mu=0$ is not possible.
\end{lemma}

\begin{proof}
Define
\begin{multline*}
V=\mathrm{Ker}(e_1:M_{-1}\to M_0)\cap
\mathrm{Ker}(e_{-1}e_2:M_{-1}\to M_0)\cap\\ \cap
\mathrm{Ker}(e_1e_{-2}e_2:M_{-1}\to M_0),
\end{multline*}
\begin{multline*}
W=\mathrm{Ker}(e_{-1}:M_{1}\to M_0)\cap
\mathrm{Ker}(e_{1}e_{-2}:M_{1}\to M_0)\cap\\ \cap
\mathrm{Ker}(e_{-1}e_{2}e_{-2}:M_{1}\to M_0).
\end{multline*}
Since $\dim M_{-1}=\infty$ and $\dim M_0<\infty$, $V$ is a vector subspace
of finite codimension in $M_{-1}$. Since $\dim M_{1}=\infty$ and $\dim M_0<\infty$,
$W$ is a vector subspace of finite codimension in $M_{1}$.  In order not to get a
direct contradiction using the Principal Observation, we assume $e_2v\neq 0$ for all
$0\neq v\in V$ and $e_{-2}w\neq 0$ for all  $0\neq w\in W$.
Then $\dim e_2V=\infty$ and, by Lemma~\ref{l2}, for every $0\neq v\in V$ we have
$(e_1^3-6e_2e_1+6e_3)e_2v=0$.

Since the codimension of $W$ in $M_1$ is finite, the intersection
$W'=e_2V\cap W$ is infinite-dimensional. Note that
\begin{displaymath}
e_{-1}(e_1^3-6e_2e_1+6e_3)=
6e_1^2e_0+6e_1^2-12e_2e_0-18e_1^2+24e_2\mod U(\V)e_{-1}.
\end{displaymath}
Choose $v\in V$ such that $w_v:=e_2v\in W'\setminus\{0\}$.  The
equality $e_{-1}(e_1^3-6e_2e_1+6e_3)e_2v=0$ implies that
$(2e_2-e_1^2)w_v=0$. In particular, for this $v$ we have
$e_{-2}(2e_2-e_1^2)w_v=0$. However,
\begin{displaymath}
e_{-2}(2e_2-e_1^2)=2e_2e_{-2}+2e_0-c-e_1^2e_{-2}-6e_1e_{-1},
\end{displaymath}
and since $e_1e_{-2}w_v=e_{-1}w_v=0$ by assumptions, we get
$e_2e_{-2}w_v=\tau w_v$ for some $\tau\in \C$. In order not to get 
a direct contradiction using the Principal Observation, we must
assume $e_{-2}w_v\ne0$. Since $e_1e_{-2}w_v=0$, we also must assume
$e_2e_{-2}w_v\ne0$, that is $\tau\ne0$.

Denote $y=w_v$ and $x=e_{-2}y$.  Let us sum up, what we know about
$x$ and $y$:
\begin{equation}\label{eq10}
e_1x=0,\quad e_{-1}y=0,\quad x=e_{-2}y,\quad \tau y=e_2x.
\end{equation}
Let $U_+$ and $U_-$ denote the subalgebras of $U(\V)$, generated
by $e_1,e_2$ and $e_{-1},e_{-2}$ respectively. Consider the vector
space
\begin{displaymath}
N=U_-x\oplus U_+y\subset M.
\end{displaymath}
From the definition it follows that both $U_+$ and $U_-$ are
stable under the adjoint action of $e_0$. Since both $x$ and $y$
are eigenvectors for $e_0$, we derive that $N$ decomposes into a
direct sum of weight spaces which are obviously
finite-dimensional. Hence, to complete the proof we have just to
show that $N$ is stable under the action of the whole $\V$. Since
$\V$ is generated by $e_1$, $e_{-1}$, $e_2$, $e_{-2}$, it is
enough to show that $N$ is stable under the action of these four
operators. Because of the symmetry of our situation, it is even
enough to show that $N$ is stable under the action of, say $e_1$
and $e_2$.

That $e_1U_+y\subset U_+y$ and $e_2U_+y\subset U_+y$ is clear. Let us show that
$e_1U_-x\subset U_-x$. For any $a\in U_-$ we have $e_1ax=ae_1x+[e_1,a]x$.
By \eqref{eq10}, $ae_1x=0$. Further, $[e_1,a]=\sum_{i,j}a_{i,j}e_0^ic^j$ for some
$a_{i,j}\in U_-$. Since $x\in M_{-1}$, we have $e_0^ic^jx=\xi x$ for some
$\xi\in\C$. Therefore $[e_1,a]x\in U_-x$, which means that $e_1$ preserves $U_-x$
and hence $N$.

Finally, let us show that $e_2U_-x\subset N$. For any $a\in U_-$
we have $e_2ax=ae_2x+[e_2,a]x$.  By \eqref{eq10}, $e_2x=\tau y\neq
0$. Let $A=e_{i_1}\dots e_{i_l}$ be a monomial, where
$i_s\in\{-1,-2\}$ for all $s=1,\dots,l$. If $i_l=-1$ we have
$Ae_2x=0$ since $e_{-1}y=0$. If $i_l=-2$ we have $Ae_2x=\zeta
e_{i_1}\dots e_{i_{l-1}}x\in U_-x$ for some $\zeta\in\C$ by
\eqref{eq10}. This implies that $ae_2x\in N$. Let us write the
element $[e_2,a]$ in the PBW basis corresponding to the order
$\dots,e_{-2},e_{-1},e_0,c,e_1$. By \eqref{eq10}, $e_1x=0$, and
hence all terms, which end on $e_1$ will vanish. This means that
$[e_2,a]x=\sum_{i,j}a_{i,j}e_0^ic^jx$ for some $a_{i,j}\in U_-$.
In the previous paragraph it was shown that in this case
$[e_2,a]x\in U_-x$. This completes the proof.
\end{proof}

\section{Corollaries from Theorem~\ref{tmain}}\label{s3}

As an immediate corollary from Theorem~\ref{tmain} we have:

\begin{corollary}\label{c2}
Let $M$ be a simple weight $\V$-module. Assume that there exists
$\lambda\in \C$ such that $0<\dim M_{\lambda}<\infty$. Then $M$ is
a Harish-Chandra module. Consequently, $M$ is either a simple
highest or lowest weight module or a simple module from the
intermediate series.
\end{corollary}

\begin{proof}
Assume that this $M$ is not a Harish-Chandra module. Then there should exists
$i\in\Z$ such that $\dim M_{\lambda+i}=\infty$. In this case
Theorem~\ref{tmain} implies $\dim M_{\lambda}=\infty$, a contradiction.
Hence $M$ is a Harish-Chandra module, and the rest of the statement follows
from \cite[Theorem~1]{M}.
\end{proof}

The following corollary gives a positive answer to \cite[Problem~3.3]{X}:

\begin{corollary}\label{c3}
Every pointed $\V$-module is a Harish-Chandra module.
\end{corollary}

\begin{proof}
Every pointed module satisfies the conditions of Corollary~\ref{c2}
by definition. Hence the statement follows from Corollary~\ref{c2}.
\end{proof}

The following corollary gives a positive answer to \cite[Conjecture~1]{Ma1}:

\begin{corollary}\label{c4}
There are no simple mixed $\V$-modules.
\end{corollary}

\begin{proof}
Let $M$ be a simple mixed $\V$-module. Then, by the definition, there exists
$\lambda\in \C$ and $i\in \Z$ such that
$\dim M_{\lambda}=\infty$ and $\dim M_{\lambda+i}<\infty$. However,
Theorem~\ref{tmain} implies $\dim M_{\lambda+i}=\infty$. A contradiction.
\end{proof}

%\vspace{1cm}
%
%\noindent
%{\bf Acknowledgments.} The research was done during the visit of the first author to
%Beijing Normal University in May 2005. The hospitality and
%financial support of Beijing Normal University are gratefully
%acknowledged. For the first author the research was partially
%supported by the Royal Swedish Academy of Sciences, the Swedish
%Research Council, and the Swedish Foundation for International
%Cooperation in research and Higher Education (STINT). The second
%author is supported by  NSERC, and NSF of China (Grant 10371120
%and 10431040).

\vspace{1cm}

\noindent
V.M.: Department of Mathematics, Uppsala University, Box 480,
SE-751 06, Uppsala, SWEDEN; e-mail: {\tt mazor\symbol{64}math.uu.se}
\vspace{0.5cm}

\noindent
K.Z.: Department of Mathematics, Wilfrid Laurier University, Waterloo,
Ontario, N2L 3C5, Canada; and Institute of Mathematics, Academy of
Mathematics and System Sciences, Chinese Academy of Sciences, Beijing
100080, PR China; e-mail: {\tt kzhao\symbol{64}wlu.ca}

\end{document}